\documentclass{amsart}
\usepackage{amsmath,amssymb,amsthm,amsfonts,amstext}

\newtheorem{theorem}{Theorem}

\newtheorem{lemma}[theorem]{Lemma}

\newtheorem{proposition}[theorem]{Proposition}

\newcommand{\cqd}{\nopagebreak\hfill\fbox{ }}

\newcommand{\LL}{\mathcal{L}}

\begin{document}

\title{Product type potential on the $X\,Y$ model: selection of maximizing probability and a large deviation principle}

\author{\hspace{1cm}J. Mohr\\
\centerline{IME, UFRGS - Porto Alegre, Brasil}
\bigskip
} 

\begin{abstract}

	Given an interval $[a,b]$  the associated $X\,Y$ model is the space $\Omega=[a,b]^\mathbb{N}$ with an a priori probability $\nu $  on the state space $[a,b]$. 



We will present here the case of the product type potential on the $X\,Y$ model and in this setting we can show the explicit expression of the equilibrium probability.

We will also consider questions about Ergodic Optimization,  maximizing probabilities, subactions and
we will show selection of a maximizing probability, when temperature goes to zero. 

Finally we show a large deviation principle  when temperature goes to zero and  we present an explicit expression for the deviation function.

\end{abstract}

\date{\today}

\maketitle

\bigskip
\thanks{ J. Mohr is partially supported by CNPq \, .}


\section{Introduction}
Let $\Omega=[0,1]^{\mathbb{N}}$ be the symbolic
space $X\,Y$ and the a priori probability $d\,a$ (Lebesgue).

We consider the metric in $\Omega=[0,1]^{\mathbb{N}}$ given by:
\[d(x,y) = \sum_{n=1}^{\infty}\frac{|x_n-y_n|}{2^{n}}\]
where $x=(x_1,x_2,...)$ and $y=(y_1,y_2,...)$ are on $\Omega$.
Note that $\Omega$ is compact by Tychonoff's theorem.
We denote by $\mathcal{C}$ the space of continuous functions from $\Omega\to \mathbb{R}$.

Given a continuous function $f:\Omega \to \mathbb{R}$
let $\mathcal{L}_{f}:\mathcal{C}\to\mathcal{C}$ be the Ruelle operator
that sends $\varphi\mapsto \mathcal{L}_{f}(\varphi)$,
which is defined for each $x\in\Omega$ by the following
expression
\begin{equation} \label{mainR}
\mathcal{L}_{f}(\varphi)(x)
=
\int_0^1 \,
e^{ f(a,x_1,x_2,...)}\,
\varphi(a,x_1,x_2,...)\,d\,a.
\end{equation}
As usual, we define the dual of the Ruelle operator, denoted by $\mathcal{L}_{f}^*$, on the space of Borel measures on $\Omega$ as the operator that send a measure $\mu$ to the measure $\mathcal{L}_{f}^*(\mu)$ defined, for each $\varphi\in\mathcal{C}$, by  
$$\int_{\Omega}\varphi\, d\,\mathcal{L}_{f}^*(\mu)= \int_{\Omega}\mathcal{L}_{f}(\varphi)\, d\,\mu.$$ 

The general case, where $\Omega=M^{\mathbb{N}} $, M is a compact set and the a priori probability is not necessarily Lebesgue is studied in \cite{LMMS} and is called one-dimensional lattice system theory. If we  suppose  $M=[0,1]$ and the a priori probability is Lebesgue, this is so-called $X\,Y$ one-dimensional model (see \cite{BCLMS}).
It is a classical problem in Physics to analyze
the Statistical Mechanics of lattices when the spin are on $S^1$ (see \cite{FH}).
It is shown in \cite{LMMS} (see Theorem 1 and 2)  that  if $f$ is Lipschitz then there exists a strictly positive Lipschitz eingenfunction $h_f$ for $\mathcal{L}_{f}$ associated to a positive eigenvalue $\lambda_f$ and also the existence of an eigenprobability  for $\mathcal{L}_{f}^*$. Moreover, the eigenvalue $\lambda_f$ is simple (which means the eigenfunction is unique up to a multiplicative constant).

We denote by  $\mathcal{M}_{\sigma}$ the set of invariant measures for the shift map, $\sigma:\Omega\to \Omega$, defined by $\sigma(x_{1},x_{2},x_{3},...)=(x_2,x_3,x_4,...).$
In \cite{LMMS} (see Definition 2) was defined the entropy $h(\mu)$ of  $\mu\in\mathcal{M}_{\sigma}$ and was proved (see Theorem 3) a variational principle: given a Lipschitz   potential $f$ and $\lambda_f$ is  the maximal eigenvalue of $\mathcal{L}_{f}$ then 
$$\log{\lambda_f}=\sup_{\mu\in\mathcal{M}_{\sigma}}\bigg\{h(\mu)+\int_{\Omega} f(x)d\mu(x)\bigg\}. $$ 
Moreover the supremum is attained on  the eigenprobability of the dual of the Ruelle operator.

These are theoretical questions on the Thermodynamic Formalism for the $X\,Y$ model which were already addressed on some recent papers. However, there is lack of interesting examples where the theory can be applied.
Here we will present several results and explicit examples on the Thermodynamic Formalism of the $X\,Y$ model in order to fill this gap.

\medskip

We consider a   continuous potential
$f:\Omega \to \mathbb{R}$  of the form
\[f(x)\,=\,
f(x_1,x_2,x_3,\ldots) = \sum_{j=1}^\infty f_j(x_j)
\]
where $f_j:[0,1]\to\mathbb{R}$ are fixed functions. We say that the  function $f$ is of the product type. We will also suppose  that
 $ \sum_{j=1}^\infty f_j(x_j)$ is absolutely convergent, for all $x\in\Omega$.
 
 We will assume in some examples that each function $f_j$, $j \in \mathbb{N}$, is a Lipschitz functions with Lipschitz constant smaller than $\frac{1}{2^j}$. In this case one can show that $f: \Omega \to \mathbb{R}$   is Lipschitz.

Functions of the product type are studied in \cite{CDLS} in the case $\Omega=M^{\mathbb{N}} $  where M is a finite or countable alphabet. In \cite{CDLS} was shown, among other things,  explicit formulae for the leading eigenvalue, the eigenfunction and eigenmeasure of the Ruelle operator.

In section \ref{explicita} we will exhibit the explicit expression of the maximal eigenvalue, of the positive eigenfunction of the Ruelle operator and of the eigenprobability of the dual of the Ruelle operator, when $M=[0,1]$.  If $f$ is Lipschitz  we know, by \cite{LMMS}, that the eigenprobability satisfies a variational principle, and hence this measure is the equilibrium probability for $f$.

\bigskip

Let $\beta=1/T$ be the inverse of the temperature $T$, if we consider the potential $\beta f$ and we denote by $\tilde \mu_{\beta}$ the eigenprobability of $\mathcal{L}_{\beta f}^*$, its well known that the limits (in the weak* topology) of $\tilde \mu_{\beta}$,  when $\beta\to \infty$, are related with the following problem: given $f:\Omega\to \mathbb{R}$ Lipschitz continuous, we want to find probabilities that maximize $\int_{\Omega}f(x)d\mu(x)$ over $\mathcal{M}_{\sigma}$. If we define  $\displaystyle m(f)=\max_{\mu\in\mathcal M_\sigma} \left\{ \int_{\Omega} f d\mu \right\}\,,$ any measure that attains the maximal value is called a maximizing measure for $f$. See \cite{LMMS} for general results in ergodic optimization theory, when $M=[0,1]$.

It is shown in \cite{LMMS}: if for some subsequence we have $\tilde \mu_{\beta_n}\rightharpoonup \mu_{\infty}$, when $n\to \infty$, then $ \mu_{\infty}$ is a maximizing measure.

One interesting question is:  $\tilde \mu_{\beta}$ converges to a maximizing measure, when $\beta\to \infty$? In the afirmative case we say we have selection of this maximizing measure.
The problem of selection and non selection of a maximizing measure was studied in several works, see \cite{ER} and \cite{CL} for examples of non selection in the case $M$ is the unitary circle.

We will  show in section \ref{principal} that we have selection of a maximizing measure in the case $f$ is of the product type and $f(a,a,a,...)$ has one or two maximum points in $[0,1]$,  also  a large deviation principle is true for this convergence.

 In \cite{LMST} was  shown a large deviation principle in the case $M=[0,1]$ and the maximizing probability is unique for a potential that depends only in two coordinates.  In the present work we do not suppose the maximizing probability is unique and the potential can depends on all coordinates.

 \section{Explicit expressions for eigenfunction and eigenprobability of functions of product type}\label{explicita}
 
Let us consider a   continuous potential of the product type
 $f:\Omega \to \mathbb{R}$ defined by
 \[f(x)\,=\,
 f(x_1,x_2,x_3,\ldots) = \sum_{j=1}^\infty f_j(x_j),
 \]
 where $f_j:[0,1]\to\mathbb{R}$ are fixed functions and  such  that
 $ \sum_{j=1}^\infty f_j(x_j)$ is absolutely convergent, for all $x\in\Omega$.
 
 Sometimes is more convenient use the following notation: $g_i(a)=e^{f_i(a)}$, then \newline $\displaystyle e^{f(x_1,x_2,x_3,...)}=e^{\sum_{j=1}^\infty f_j(x_j)}=\prod_{j=1}^{\infty}g_j(x_j):=g(x_1,x_2,x_3,...).$

 In this way, equation \eqref{mainR} became
 $$\mathcal{L}_{f}(\varphi)(x)
 =
 \int_0^1 \,
 g_1(a)\prod_{j=2}^{\infty}g_j(x_{j-1})\,
 \varphi(a,x_1,x_2,...)da.$$
 In this section we will show  the explicit  expressions for the maximal eigenvalue and for positive eigenfunction of $\mathcal{L}_{f}$ and for  the eigenprobability of $\mathcal{L}^*_{f}$.
\medskip

The following proposition is the analogous of Theorem 4.1 in \cite{CDLS}.
\begin{proposition} \label{pro1} Suppose $f$ satisfies 
 $ \displaystyle \sum_j \sum_{i>j} f_i(x_j)<\infty$, for all $x=(x_1,x_2,...)\in\Omega$.
If we define  $\displaystyle h_f(x)=\prod_{j=1}^{\infty}h_j(x_j)$, where $\displaystyle h_j(b)=\prod_{i>j}g_i(b)$ and $\displaystyle\lambda_f=\int_0^1 \,\prod_{j=1}^{\infty}g_j(b)\,db$. Then $\mathcal{L}_{f}(h_f)=\lambda_f h_f$.
	
\end{proposition}

\bigskip

\textbf{Proof:}
First we will show that $\lambda_f<\infty$ and $h_f(x)<\infty$, for all $x\in \Omega$.

In fact, as  $\prod_{j=1}^{\infty}g_j(b)=e^{\sum_{j=1}^\infty f_j(b)}=e^{f(b,b,...)}$ and $b\to f(b,b,...)$ is continuous, we have that  $\lambda_f=\int_0^1 \, \prod_{j=1}^{\infty}g_j(b)\,db<\infty$.
Note that 
$\sum_j \log  h_j (x_j)=\sum_j \sum_{i>j} \log g_i(x_j)=\sum_j \sum_{i>j} f_i(x_j) <\infty$,  this implies
$\displaystyle h_f(x)=\prod_{j=1}^{\infty}h_j(x_j)<\infty$.

Now we will show that $h_f$ is a eigenfunction to  $\mathcal{L}_{f}$: as $h_j(b)=\prod_{i>j}g_i(b)  $, multiplying it by $g_j(b)$  we obtain
\begin{equation}\label{rel_gh}
g_j(b)\,h_j(b)=g_j(b)\prod_{i>j}g_i(b)=\prod_{i>j-1}g_i(b)=h_{j-1}(b).
\end{equation}
In particular $g_1(a)h_1(a)=\prod_{i=1}^{\infty}g_i(a)$ and $ g_j(x_{j-1})\,
h_j(x_{j-1})=h_{j-1}(x_{j-1})$. Hence

$$\mathcal{L}_{f}(h_f)(x)
=
\int_0^1 \,
g_1(a)\prod_{j=2}^{\infty}g_j(x_{j-1})\,
h_1(a)\prod_{j=2}^{\infty}h_j(x_{j-1})da=$$
$$\int_0^1 \,
g_1(a)h_1(a)\prod_{j=2}^{\infty}g_j(x_{j-1})\,
h_j(x_{j-1})da=\int_0^1 \,\prod_{i=1}^{\infty}g_i(a)da
\prod_{j=2}^{\infty}\,
h_{j-1}(x_{j-1})=\lambda_f h_f(x). $$
\cqd

\textbf{Remark}: We need not suppose that $f$ is Lipschitz to prove the previous theorem, but if we suppose $f$ Lipschitz then we have that $h_f$ is the unique strictly positive eigenfunction for $\mathcal{L}_f$ (see Theorem 1 in \cite{LMMS}).

Note that $\displaystyle h_f(x)=\prod_{j=1}^{\infty}\prod_{i>j}g_i(x_j)=\prod_{j=1}^{\infty}\prod_{i>j}e^{f_i(x_j)}=e^{\sum_{j=1}^{\infty}\sum_{i>j}f_i(x_j)}$.

\begin{proposition}
	Assume that the functions $f_j$, $j \in \mathbb{N}$ are Lipschitz with Lipschitz constant smaller than $\frac{1}{2^j}$ and that for some $\bar x$ we know that $\sum_{j=1}^{\infty}\sum_{i>j}f_i(\bar x_j)<\infty$. Then, the hypothesis of Proposition \ref{pro1} are true.
\end{proposition}

{\bf Proof:} $$\bigg|\sum_{j=1}^{\infty}\sum_{i>j}f_i( x_j)-\sum_{j=1}^{\infty}\sum_{i>j}f_i(\bar x_j)\bigg|\leq \sum_{j=1}^{\infty}\sum_{i=j+1}^{\infty}\frac{1}{2^i}=\sum_{j=1}^{\infty}\frac{1}{2^{j-1}}=2.$$
\cqd

\bigskip
We say that a potential $\tilde f$ is normalized if $\mathcal{L}_{\tilde f}(1)=1$ or $\int_0^1 \,
e^{ \tilde f(a,x_1,x_2,...)}\,
d\,a=1$. \newline
Given a potential $f$, let $h_f$ and $\lambda_f$ be  as in the Proposition \ref{pro1}, we define the normalized potential associated to $f$, as usual by $\tilde f=f+\log h_f-\log h_f\circ \sigma-\log\lambda_f$.

In the exponential scale $e^{\tilde f}$ became
 $$\tilde g(x)=\frac{g(x) h_f(x)}{h_f\circ \sigma(x) \lambda_f}=\frac{g(x) }{ \lambda_f}\prod_{i=1}^{\infty}\frac{h_i(x_i)}{h_i(x_{i+1})}=\frac{g(x)h_1(x_1) }{ \lambda_f}\prod_{i=1}^{\infty}\frac{h_{i+1}(x_{i+1})}{h_i(x_{i+1})},$$
 now using equation \eqref{rel_gh} and the definition of $g$ we get
$$\tilde g(x)=\frac{\prod_{j=1}^{\infty}g_{j}(x_{j})h_1(x_1) }{ \lambda_f \prod_{i=1}^{\infty}g_{i+1}(x_{i+1})}=\frac{g_{1}(x_{1})h_1(x_1)} { \lambda_f }=\frac{g_{1}(x_{1})g_{2}(x_{1})h_2(x_1)} { \lambda_f }=...$$
$$...=\frac{\prod_{i=1}^{\infty}g_{i}(x_{1})}{\int_0^1 \prod_{j=1}^{\infty}g_j(b)db}=\frac{e^{\sum_{i=1}^{\infty}f_{i}(x_{1})}}{\int_0^1 e^{\sum_{j=1}^{\infty}f_j(b)}db}=\frac{e^{f(x_1,x_1,x_1...)}}{\int_0^1 e^{f(b,b,b,...)}db}.$$
This implies $\tilde g$ (and $\tilde f$) depends only on the first coordinate of $x$.

\bigskip

 It is known   from \cite{LMMS} (see Theorem 2) that if $f$ is Lipschitz continuous then there exists a unique eigenprobability $\tilde\mu_f$ for $\mathcal{L}^{*}_{\tilde f}$, and that the measure $\mu_f=\frac{1}{h_f}\tilde \mu_f$ is an eingenmeasure for $\mathcal{L}^{*}_{ f}$, where $h_f$ is the unique eigenfunction of $\mathcal{L}_{ f}$ associated to the maximal eigenvalue $\lambda_f$. The next proposition exhibits the explicit form of these measures.

\medskip
The following Proposition is analogous to Theorem 4.2 of \cite{CDLS}, although the expression of $\mu_f$ is slightly different and the proof is more direct.
\begin{proposition}\label{medidas}
	Suppose $f$  satisfies the hypothesis of Proposition \ref{pro1}. Let  $\mu_f=\otimes_{n=1}^{\infty} \mu_n$ and $\tilde \mu_f=\otimes_{n=1}^{\infty}\tilde \mu_0$ be measures  of the product type given by the following expressions
	$$d\mu_n(a)=\frac{\prod_{i=1}^{n}g_{i}(a)\,d a}{\int_0^1 \prod_{j=1}^{\infty}g_j(b)db}\,\,\,, \,\,\, d\tilde\mu_0(a)=\frac{\prod_{i=1}^{\infty}g_{i}(a)\,da}{\int_0^1 \prod_{j=1}^{\infty}g_j(b)db}=\frac{e^{F(a)}\,da}{\int_0^1 e^{F(b)}db}= \tilde g(a)da,$$
where $F:[0,1]\to\mathbb{R}$ is defined by $F(a)=\sum_{i=1}^{\infty}f_i(a)=f(a,a,a,...)$.

Then, we have $\mathcal{L}^{*}_{ f}( \mu_f)=\lambda_f \mu_f$, $\mathcal{L}^{*}_{\tilde f}(\tilde \mu_f)=\tilde \mu_f$ and $\mu_f=\frac{1}{h_f}\tilde \mu_f$.
\end{proposition}

\textbf{{Proof:}}
Note that, by definition, $d\mu_1(a)=\frac{g_1(a)da}{\lambda_f}$ and $d\mu_{n+1}(x_n)=\frac{\prod_{i=1}^{n+1}g_{i}(x_n)\,d x_n}{\int_0^1 \prod_{j=1}^{\infty}g_j(b)db}=g_{n+1}(x_n)d\mu_n(x_n)$, we will use these equalities in the following calculation.

Let $\varphi:[0,1]^{\mathbb{N}}\to\mathbb{R}$, then
$$\int_{[0,1]^{\mathbb{N}}} \varphi\, d\mathcal{L}^{*}_{ f}( \mu_f)=\int_{[0,1]^{\mathbb{N}}}\mathcal{L}_{ f}(\varphi )  d\mu_f=$$
$$=\int_{[0,1]^{\mathbb{N}}}\int_{[0,1]}\varphi(a,x_1,x_2,...)g_1(a)g_2(x_1)g_3(x_2)...da\, d\mu_1(x_1)\,d\mu_2(x_2)\,...= $$
$$=\int_{[0,1]^{\mathbb{N}}}\int_{[0,1]}\varphi(a,x_1,x_2,...)g_1(a)da\, g_2(x_1)d\mu_1(x_1)\,g_3(x_2)d\mu_2(x_2)\,...= $$
$$=\lambda_f\int_{[0,1]^{\mathbb{N}}}\int_{[0,1]}\varphi(a,x_1,x_2,...)\frac{g_1(a)da}{\lambda_f}\, g_2(x_1)d\mu_1(x_1)\,g_3(x_2)d\mu_2(x_2)\,...= $$

$$=\lambda_f\int_{[0,1]^{\mathbb{N}}}\varphi(a,x_1,x_2,...)d\mu_1(a)\, d\mu_2(x_1)\,d\mu_3(x_2)\,...=\lambda_f \int_{[0,1]^{\mathbb{N}}} \varphi \,d\mu_f, $$
and this implies that $\mathcal{L}^{*}_{ f}( \mu_f)=\lambda_f \mu_f$.

And, as $\tilde g(a)d a=d\tilde \mu_0(a)$ we have
$$\int_{[0,1]^{\mathbb{N}}} \varphi\, d\mathcal{L}^{*}_{\tilde f}(\tilde \mu_f)=\int_{[0,1]^{\mathbb{N}}}\mathcal{L}_{\tilde f}(\varphi )  d\tilde\mu_f=$$
$$=\int_{[0,1]^{\mathbb{N}}}\int_{[0,1]}\tilde g(a)\varphi(a,x_1,x_2,...)da\, d\tilde\mu_0(x_1)\,d\tilde\mu_0(x_2)\,...= $$
$$=\int_{[0,1]^{\mathbb{N}}}\int_{[0,1]}\varphi(a,x_1,x_2,...)\tilde g(a)da\, d\tilde\mu_0(x_1)\,d\tilde\mu_0(x_2)\,...=$$
$$=\int_{[0,1]^{\mathbb{N}}}\varphi(a,x_1,x_2,...)d\tilde\mu_0(a)\, d\tilde\mu_0(x_1)\,d\tilde\mu_0(x_2)\,...=\int_{[0,1]^{\mathbb{N}}} \varphi\,  d \tilde u_f,$$
this implies that $\mathcal{L}^{*}_{\tilde f}(\tilde \mu_f)=\tilde \mu_f$.

Finally, as  $\frac{1}{h_j(x_j)}\tilde g(x_j)dx_j=  \frac{1}{\prod_{i>j}^{\infty}g_{i}(x_j)}\,\frac{\prod_{i=1}^{\infty}g_{i}(x_j)dx_j}{\int_0^1 \prod_{j=1}^{\infty}g_j(b)db}     =\frac{\prod_{i=1}^{j}g_{i}(x_j)dx_j}{\int_0^1 \prod_{j=1}^{\infty}g_j(b)db}=d\mu_j(x_j)$, 
we get $\mu_f=\frac{1}{h_f}\tilde \mu_f$, because 
 $$\int_{\Omega} \varphi(x)\frac{1}{h_f}d\tilde \mu_f(x)=\int_{\Omega} \varphi(x_1,x_2,...)\frac{1}{h_1(x_1)}\tilde g(x_1)dx_1\frac{1}{h_2(x_2)}\tilde g(x_2)dx_2...=$$
 $$=\int_{\Omega} \varphi(x_1,x_2,...)d\mu_1(x_1)d\mu_2(x_2)...=\int_{\Omega}\varphi(x)d\mu_f(x).$$
\cqd

\medskip

\textbf{Remark}: We need not suppose that $f$ is Lipschitz to prove the previous theorem, but if we suppose $f$ Lipschitz then we have that $\tilde \mu_f$ is the unique fixed point to $\mathcal{L}^*_{\tilde f}$ (see Theorem 2 in \cite{LMMS}). The hypothesis of Proposition \ref{pro1} is used, in the previous theorem, only to prove that $\mu_f=\frac{1}{h_f}\tilde \mu_f$.


  Following  \cite{LMMS} (see Definition 1), if $\tilde f$ is Lipschitz continuous and normalized, and $\tilde \mu_f$ is such that $\mathcal{L}^{*}_{\tilde f}(\tilde \mu_f)=\tilde \mu_f$,  the entropy of $\tilde\mu_f$ is definided  by,
$$ h(\tilde\mu_f)=-\int_{[0,1]^{\mathbb{N}}}\tilde f(x)d\tilde \mu_f(x) =$$
$$=-\int_{[0,1]^{\mathbb{N}}}\log \tilde g(x)d\tilde \mu_f(x)=-\int_{[0,1]^{\mathbb{N}}}\log  \frac{\prod_{i=1}^{\infty}g_i(x_1)}{\int_0^1\prod_{i=1}^{\infty}g_i(b)db}  d\tilde \mu_f(x)=$$
$$=-\int_{[0,1]^{\mathbb{N}}}\bigg[\log  {\prod_{i=1}^{\infty}g_i(x_1)}-\log {\int_0^1\prod_{i=1}^{\infty}g_i(b)db}\bigg]  \tilde g(x_1)dx_1  \tilde g(x_2)dx_2 ...\tilde g(x_n)dx_n...=$$
$$=\log \lambda_f-\int_{[0,1]}\log  {\prod_{i=1}^{\infty}g_i(x_1)}  \tilde g(x_1)dx_1=\log\lambda_f-\int_{[0,1]} \sum_{i=1}^{\infty} f_i(a) \tilde g(a)da. 
$$
This is an explicit expression for the entropy of this example.

\bigskip 

Also we compute \newline
	$$\int_{[0,1]^{\mathbb{N}}} f d\tilde \mu_f=\int_{[0,1]^{\mathbb{N}}} \sum_{i=1}^{\infty} f_i(x_i) \tilde g(x_1)dx_1   ...\tilde g(x_i)dx_i...=$$
	$$=\sum_{i=1}^{\infty}\int_{[0,1]^{\mathbb{N}}} f_i(x_i) \tilde g(x_1)dx_1 ...\tilde g(x_i)dx_i...=\sum_{i=1}^{\infty}\int_{[0,1]} f_i(x_i) \tilde g(x_i)dx_i=\int_{[0,1]} \sum_{i=1}^{\infty}f_i(a) \tilde g(a)da.$$
\bigskip

And this implies that $ h(\tilde\mu_f)=\log\lambda_f-\int_{[0,1]^{\mathbb{N}}} f d\tilde \mu_f  $ or  $$\log\lambda_f= h(\tilde\mu_f)+\int_{\Omega} f d\tilde \mu_f . $$

This shows  that $\tilde \mu_f$ satisfies a variational principle, as in Theorem 3 of  \cite{LMMS}, i.e., let $f$ be a Lipschitz continuous potential and $\lambda_f$ be  the maximal eigenvalue of $\LL_f$, then
$$  \log \lambda_f\,=\,P(f)\,=\,\sup_{\mu\in\mathcal{M}_{\sigma}} \bigg\{h(\mu)+\int_{\Omega} f (x) d\mu(x)\bigg\}.  $$
And the supremum is attained on the  measure $\tilde \mu_f$.

\bigskip
\section{Zero temperature, selection of the maximizing measure and large deviation principle}\label{principal}
Now we will analyze the question  of zero temperature, when $\beta\to \infty$, for this example. General results on Ergodic Optimization and selection when temperature goes to zero, for the case $\Omega=\{1,...,d\}^{\mathbb{N}}$, can be found in \cite{BLL}.

 For each $\beta>0$  we consider the potential  $\beta f(x)=\sum_{j=1}^\infty\beta f_j(x_j)$, where $\beta f$ is Lipschitz and satisfies the hypothesis of Proposition \ref{pro1}, so the eigenfunction of $\mathcal{L}_{\beta f}$ is given by  $h_{\beta}(x)=e^{\sum_{j=1}^{\infty}\sum_{i>j}\beta f_i(x_j)}=e^{\beta\sum_{j=1}^{\infty}\sum_{i>j} f_i(x_j)} $.
And, the equilibrium probability is given by $\tilde\mu_{\beta}=\otimes_{n=1}^{\infty} \tilde\mu_{0,\beta}$ where, by Proposition \ref{medidas},
\begin{equation}\label{mu0}
d\tilde\mu_{0,\beta}(a)=\frac{ e^{\beta\sum_{i=1}^{\infty}f_i(a)}}{\int_0^1 e^{\beta\sum_{i=1}^{\infty} f_i(b)}db}\,da=\frac{ e^{\beta F(a)}}{\int_0^1 e^{\beta F(b)}db}\,da.
\end{equation}

As usual, we would like to  investigate the limits of $\tilde\mu_{\beta}$ and $ \frac{1}{\beta} \log h_{\beta}(x)$, when $\beta \to \infty $.

The limits of $\tilde \mu_{\beta }$ are related  with the following problem: for $f:\Omega \to \mathbb{R}$ fixed above, we want to find probabilities that  maximize the value
$ \displaystyle\int_{\Omega} f(x) \,d \mu( {x}). $

 We define  
$$m(f)=\max_{\mu\in\mathcal M_\sigma} \left\{ \int_{\Omega} f d\mu \right\}\,.$$

Any of the probability measures which attains the maximal value will be called a maximizing probability measure, which will be  denoted generically by $\mu_\infty$.
\bigskip

We say that $u$ is a calibrated subaction if
\begin{equation}\label{subaction}
m(f)=\max_{a\in[0,1]}\{f(ax)+u(ax)-u(x)\}.
\end{equation}

We know that by Proposition 10 in \cite{LMMS} that, if the potential $f$ is Lipschitz, then\newline  i)$\hspace{4cm}\displaystyle\lim_{\beta \to \infty}  \frac{1}{\beta}\log\lambda_{\beta}= m(f),$

 where $\lambda_{\beta}=\int_0^1 e^{\beta\sum_{j=1}^\infty  f_j(a)}da$.\newline
 ii) Any limit, in the uniform topology,
 $$u:=\lim_{n\to\infty}\frac{1}{\beta_n}\log(h_{\beta_n }),$$
 is  a calibrated subaction for $f$.

\medskip

Note that
$\displaystyle\frac{1}{\beta} \log h_{\beta}(x)=\sum_{j=1}^{\infty}\sum_{i>j} f_i(x_j)$
does not depends on $\beta$, hence by the previous result we have that $u(x)=\sum_{j=1}^{\infty}\sum_{i>j} f_i(x_j)$ is a calibrated subaction. 

\bigskip
\begin{proposition}
	$\displaystyle m(f)=\max_{a\in[0,1]}\sum_{i=1}^{\infty}f_i(a).  $
\end{proposition}
\textbf{{Proof:}}
 Let $u(x)=\sum_{j=1}^{\infty}\sum_{i>j} f_i(x_j)$ be a calibrated subaction,  first see that:\newline  $u(x)=\sum_{i=2}^{\infty}f_i(x_1)+\sum_{i=3}^{\infty}f_i(x_2)+\sum_{i=4}^{\infty}f_i(x_3)+...$, and \newline $u(ax)=\sum_{i=2}^{\infty}f_i(a)+\sum_{i=3}^{\infty}f_i(x_1)+\sum_{i=4}^{\infty}f_i(x_2)+...$, hence \newline
$u(ax)-u(x)=\sum_{i=2}^{\infty}f_i(a)+\sum_{i=3}^{\infty}f_i(x_1)-\sum_{i=2}^{\infty}f_i(x_1)+
\sum_{i=4}^{\infty}f_i(x_2)-\sum_{i=3}^{\infty}f_i(x_2)+...=\sum_{i=2}^{\infty}f_i(a)-\sum_{i=2}^{\infty}f_i(x_{i-1}).$	Therefore
$$ f(ax)+u(ax)-u(x)= f_1(a)+ \sum_{i=2}^{\infty}f_i(x_{i-1}) +\sum_{i=2}^{\infty}f_i(a)-\sum_{i=2}^{\infty}f_i(x_{i-1})= \sum_{i=1}^{\infty}f_i(a).$$
Finnaly, using equation \eqref{subaction}, we have $$m(f)=\max_{a\in[0,1]}\{f(ax)+u(ax)-u(x)\}=\max_{a\in[0,1]}\sum_{i=1}^{\infty}f_i(a).$$
	\cqd

	%


\bigskip
\begin{lemma}\label{assintotico}
	Suppose $l(\beta) =\int_{\alpha}^{\delta} e^{\beta F(t)dt}$, where $\beta$ is real and positive, $F(t), F'(t)$ and $F''(t)$ are real and continuous in $\alpha\leq t\leq \delta$. Let $t=a$ be the only point of maximum of $F(t)$ in $[\alpha,\delta]$, with $\alpha<a<\delta$, thus the asymptotic approximation as $\beta\to\infty$ is
	$$\int_{\alpha}^{\delta} e^{\beta F(t)}dt=e^{\beta F(a)}\bigg[\bigg(\frac{-2\pi}{\beta F''(a)}\bigg)^{\frac{1}{2}}+ O(\beta^{-\frac{3}{2}})\bigg].$$
\end{lemma}

\medskip

For the proof of  this Lemma see section 2.2 of \cite{Mu}.

\medskip

We will use   Lemma \ref{assintotico} to show that we have selection of the maximizing measure in the following cases:

\begin{theorem}
Let $F(b)=f(b,b,b,...)=\sum_{i=1}^{\infty}f_i(b)$ and suppose $F(b), F'(b)$ and $F''(b)$ are real and continuous in $0\leq b\leq 1$.

a) Suppose  $F$ has only one maximum in $a_1\in(0,1)$ then   $\lim_{\beta\to \infty}\tilde\mu_{0,\beta}=\delta_{a_1}$ and $\lim_{\beta\to \infty}\tilde\mu_{\beta}=\otimes_{n=1}^{\infty} \delta_{a_1}$.

b) Suppose $F$ has  two maximum   points  in $(0,1)$, say $0<a_1<a_2<1$, then we have $\lim_{\beta\to \infty}\tilde\mu_{0,\beta}= \tilde\mu_{0,\infty}=p_1\delta_{a_1}+p_2\delta_{a_2}$ and $\lim_{\beta\to \infty}\tilde\mu_{\beta}=\otimes_{n=1}^{\infty} p_1\delta_{a_1}+p_2\delta_{a_2}$, where  $p_1+p_2=1$ and $\frac{p_1}{p_2}=\sqrt{\frac{F''(a_2)}{F''(a_1)}}$.
\end{theorem}

\textbf{{Proof:}}

a) As $\tilde\mu_{\beta}=\otimes_{n=1}^{\infty} \tilde\mu_{0,\beta}$, using equation \eqref{mu0}, we need to analyse the limit of \newline $\tilde\mu_{0,\beta}(da)=\frac{ e^{\beta F(a)}}{\int_0^1 e^{\beta F(b)}db}\,da.$

If $F$ has only one maximum in $a_1\in(0,1)$, by Lemma \ref{assintotico}
	$$\int_{0}^{1} e^{\beta F(b)}db=e^{\beta F(a_1)}\bigg[\bigg(\frac{-2\pi}{\beta F''(a_1)}\bigg)^{\frac{1}{2}}+ O(\beta^{-\frac{3}{2}})\bigg].$$

Therefore, for each $a\in [0,1]$ we have
$$ \frac{ e^{\beta F(a)}}{\int_0^1 e^{\beta F(b)}db}=\frac{ e^{\beta F(a)}}{e^{\beta F(a_1)}\bigg[\bigg(\frac{-2\pi}{\beta F''(a_1)}\bigg)^{\frac{1}{2}}+ O(\beta^{-\frac{3}{2}})\bigg]}=\frac{ e^{\beta( F(a)-F(a_1))}}{\bigg(\frac{-2\pi}{\beta F''(a_1)}\bigg)^{\frac{1}{2}}+ O(\beta^{-\frac{3}{2}})} . $$

We conclude that the above expression goes to 0 if $a\neq a_1$ and goes to $\infty$ if $a=a_1$, when $\beta\to\infty$. Hence, $\lim_{\beta\to \infty}\tilde\mu_{0,\beta}(a)=\delta_{a_1}$ and $\lim_{\beta\to \infty}\tilde\mu_{\beta}=\otimes_{n=1}^{\infty} \delta_{a_1}$ and we have selection of the maximizing measure.

b) Now we consider the case where $F$ has two maximum   points  in $(0,1)$, say $0<a_1<a_2<1$, we  divide $[0,1]$ in two intervals, each one containing only one maximum point, applying the Lemma \ref{assintotico} in each interval, we obtain, as $F(a_1)=F(a_2),$
	$$\int_{0}^{1} e^{\beta F(b)}db=e^{\beta F(a_1)}\bigg[\bigg(\frac{-2\pi}{\beta F''(a_1)}\bigg)^{\frac{1}{2}}+ O(\beta^{-\frac{3}{2}})\bigg]+e^{\beta F(a_2)}\bigg[\bigg(\frac{-2\pi}{\beta F''(a_2)}\bigg)^{\frac{1}{2}}+ O(\beta^{-\frac{3}{2}})\bigg]=$$
$$=e^{\beta F(a_1)}\bigg[\bigg(\frac{-2\pi}{\beta F''(a_1)}\bigg)^{\frac{1}{2}}+\bigg(\frac{-2\pi}{\beta F''(a_2)}\bigg)^{\frac{1}{2}}+ O(\beta^{-\frac{3}{2}})\bigg], $$
hence for each $a\in[0,1]$ we have
$$ \frac{ e^{\beta F(a)}}{\int_0^1 e^{\beta F(b)}db}=\frac{ e^{\beta( F(a)-F(a_1))}}{\bigg(\frac{-2\pi}{\beta F''(a_1)}\bigg)^{\frac{1}{2}}+\bigg(\frac{-2\pi}{\beta F''(a_2)}\bigg)^{\frac{1}{2}}+ O(\beta^{-\frac{3}{2}})}.  $$
Therefore, if $a\neq a_1$ and $a\neq a_2$ the density of $\tilde \mu_{0,\beta}$ goes to 0.

 Let us fix  $\varepsilon_1,\varepsilon_2>0$  such that $a_2\notin(a_1-\varepsilon_1,a_1+\varepsilon_1)$ and $a_1\notin(a_2-\varepsilon_2,a_2+\varepsilon_2)$, we have, for  $i=1,2$, that $$\tilde\mu_{0,\beta}(a_i-\varepsilon_i,a_i+\varepsilon_i)=\frac{\int_{a_i-\varepsilon_i}^{a_i+\varepsilon_i} e^{\beta F(a)}da}{\int_0^1 e^{\beta F(b)}db}\,. $$ Now we apply Lemma \ref{assintotico} for each interval $I_i=(a_i-\varepsilon_i,a_i+\varepsilon_i)$, $i=1,2$, to obtain
$$\frac{\tilde\mu_{0,\beta}(I_1)}{\tilde\mu_{0,\beta}(I_2)}=\frac{\int_{a_1-\varepsilon_1}^{a_1+\varepsilon_1} e^{\beta F(a)}da}{\int_{a_2-\varepsilon_2}^{a_2+\varepsilon_2} e^{\beta F(a)}da}=\frac{e^{\beta F(a_1)}\bigg[\bigg(\frac{-2\pi}{\beta F''(a_1)}\bigg)^{\frac{1}{2}}+ O(\beta^{-\frac{3}{2}})\bigg]}{e^{\beta F(a_2)}\bigg[\bigg(\frac{-2\pi}{\beta F''(a_2)}\bigg)^{\frac{1}{2}}+ O(\beta^{-\frac{3}{2}})\bigg]}\approx \sqrt{\frac{F''(a_2)}{F''(a_1)}}.$$

This implies that $\tilde\mu_{0,\beta}\rightharpoonup \tilde\mu_{0,\infty}=p_1\delta_{a_1}+p_2\delta_{a_2}$, where $p_1+p_2=1$ and $\frac{p_1}{p_2}=\sqrt{\frac{F''(a_2)}{F''(a_1)}}$. And therefore we have selection of the maximizing measure.
\cqd

\bigskip
We can also prove a large deviation principle and exhibit the deviation function:

\begin{proposition} We denote $ u(x)=\sum_{j=1}^{\infty}\sum_{i>j} f_i(x_j)$  the calibrated subaction, where $x=(x_1,x_2,,...,x_j,...)$. Consider the function
	$$I(x)=\sum_{j\geq 1}u(\sigma^j(x))-u(\sigma^{j-1}(x))-f(\sigma^{j-1}(x))+m(f), $$ then
	\newline
i)$ \,\displaystyle I(x)=\sum_{j=1}^{\infty}\Big(\sum_{i=1}^{\infty}- f_i(x_j)+m(f)\Big)=\sum_{j=1}^{\infty}\Big(- F(x_j)+m(f)\Big).$
\newline
ii) For each cylinder $D=A_1\times ...\times A_n$, where $A_i$ are intervals of $[0,1]$, the following limit exists
$$\lim_{\beta\to\infty}\frac{1}{\beta}\log \tilde\mu_{\beta}(D)=-\inf_{x\in D}I(x).$$ 	
\end{proposition}

\textbf{Proof: }
i) Follows by a straigth forward calculation using the definition of $u$.



 \noindent ii) Let $D=A_1\times ...\times A_n$ be a cylinder of $\mathbb{R}^n$, then 
$\tilde\mu_{\beta}(D)=\tilde\mu_{0,\beta}(A_1)...\tilde\mu_{0,\beta}(A_n)$, using the second equality of equation \eqref{mu0} for each $\tilde\mu_{0,\beta}(A_i)$, $i=1,2,...,n$, a straight forward calculation, using that $\displaystyle m(f)=\max_{a\in[0,1]}F(a)$,  shows that
$$\lim_{\beta\to\infty}\frac{1}{\beta}\log \tilde\mu_{\beta}(D)=-\inf_{x_1\in A_1,..., x_n\in  A_n}\sum_{j=1}^{n}\Big(-F(x_j)+\max_{a\in[0,1]}F(a)\Big)=-\inf_{x\in D}I(x).$$



\cqd
\bigskip


Note that as  $\displaystyle I(x)=\sum_{j=1}^{\infty}\Big(-F(x_j)+\max_{a\in[0,1]}F(a)\Big),$
this implies that $I(x)\geq 0$, and $I(x_1,x_2,...,x_j,...)=0$, if and only if,  $x_j\in\mbox{argmax}\;F$, for all $j\in\mathbb{N}$.

Note that $I(x_1,...,x_n,x_1,...,x_n,x_1,...,x_n,...)=\infty$, if there exists $x_j\notin\mbox{argmax}\;F$, $1\leq j\leq n$.

Note also that to have $I(x)< \infty$  is necessary that $F(x_j)\to m(f)$.

\bigskip

\textbf{Example 1:} Let us define $f(x)=\sum_{i=1}^{\infty} -(x_i)^{2i}$ and suppose that we take $[-\frac{1}{2}, \frac{1}{2}]$ instead $[0,1]$. Then $f_i(a)=-a^{2i}$ and 
note that $\frac{d f_i(a)}{da}  =- 2\, i a^{2 \, i -1}$, if we define  $c_i:=\displaystyle\sup_{a\in[-\frac{1}{2}, \frac{1}{2}]}
\bigg|\frac{d f_i(a)}{da}\bigg|=2\, i 2^{-2 \, i +1}=i2^{-2 \, i +2}$, hence
$c_i< 2^{-i} $ for each $i\geq 5$. Note also that $c_i \leq 4\cdot 2^{-i}$ for $i=1,2,3$ and $4$.

Hence, we get that  the Lipschitz constant of $f_i$ is smaller than $4\cdot2^{-i}$, for all $i$.

Then,
$|f(x)-f(y)|\leq\sum_{i=1}^{\infty}|f_i(x_i)-f_i(y_i) | \leq\sum_{i=1}^{\infty}4\frac{|x_i-y_i |}{2^i} =4 d(x,y)$, i.e., $f $ is Lipschitz with constant 4.
\medskip

Also, $$F(a)=\sum_{i=1}^{\infty} f_i(a)=-\frac{1}{1-a^2}+1=1+\frac{1}{a^2-1}.$$
In this case $m(f)=F(0)=0$ and if $x=(x_1,x_2,...,x_j,....)$ then we get
$$I(x)=-\sum_{j=1}^{\infty}\Big(1+\frac{1}{x_j^2-1}\Big).$$
and
$ u(x)=\sum_{j=1}^{\infty}\sum_{i>j} f_i(x_j)=\sum_{j=1}^{\infty}\sum_{i>j} -(x_j)^{2i}$.

\bigskip

\textbf{Example 2:} Suppose we take $[-1, 1]$ instead $[0,1]$ and $f_i(a)=a^{i}i^{-\gamma},\gamma>1$,
then $$F(a)=\sum_{i=1}^{\infty} f_i(a)=\sum_{i=1}^{\infty}\frac{a^{i}}{i^{\gamma}},$$
this is the polylogarithm function.

Each $f_i$ is a Lipschitz function:
in the same way as  before we consider

$c_i:=\displaystyle\sup_{a\in[-1, 1]}
\bigg|\frac{d f_i(a)}{da}\bigg|=\sup_{a\in[-1, 1]}
\bigg|\frac{a^{i-1}}{i^{\gamma-1}}\bigg|=i^{1-\gamma}.$


The function $f$ is not Lipschitz but satisfies the hypothesis of Proposition \ref{pro1}, when $\gamma>2$. Indeed, if $\gamma>1$
$$\sum_{i>j}\frac{(x_j)^{i}}{i^{\gamma}} \leq \sum_{i>j}\frac{|(x_j)^{i}|}{i^{\gamma}} \leq \sum_{i=j+1}^{\infty}\frac{1}{i^{\gamma}} \leq \int_j^{\infty}x^{-\gamma}dx=\lim_{b\to \infty}\frac{x^{-\gamma+1}}{-\gamma+1}\bigg|_{j}^b=\frac{j^{-\gamma+1}}{-\gamma+1}  .$$ 
Now, if $\gamma>2$
$$\sum_{j=1}^{\infty}\sum_{i>j}\frac{(x_j)^{i}}{i^{\gamma}}\leq \sum_{j=1}^{\infty}\frac{j^{-\gamma+1}}{-\gamma+1} =\frac{1}{1-\gamma}+\sum_{j=2}^{\infty}\frac{j^{-\gamma+1}}{-\gamma+1}\leq\frac{1}{1-\gamma}+\int_1^{\infty}\frac{x^{-\gamma+1}}{-\gamma+1}dx=$$
$$= \frac{1}{1-\gamma}+\lim_{b\to \infty}\frac{x^{-\gamma+2}}{(-\gamma+1)(-\gamma+2)}\bigg|_{1}^b=\frac{1}{1-\gamma}+\frac{1}{(1-\gamma)(2-\gamma)}< \infty.$$ 


Note that $\max_{a\in[0,1]}F(a)$ occurs when $a=1$, hence $\max_{a\in[0,1]}F(a)=\zeta(\gamma)$.

Then,
$$I(x)=\sum_{j=1}^{\infty}\bigg(\sum_{i=1}^{\infty}- f_i(x_j)+\max_{a\in[0,1]}\sum_{i=1}^{\infty} f_i(a)\bigg)=\sum_{j=1}^{\infty}\Big(-\sum_{i=1}^{\infty}\frac{x_j^{i}}{i^{\gamma}}+\zeta(\gamma)\Big)$$
and
$$ u(x)=\sum_{j=1}^{\infty}\sum_{i>j} f_i(x_j)=\sum_{j=1}^{\infty}\sum_{i>j}\frac{x_j^{i}}{i^{\gamma}}.$$


\end{document}